\documentclass[10pt]{article}
\usepackage[utf8]{inputenc}
\usepackage[english]{babel}
\usepackage{amsmath,amsthm,amsfonts}
\usepackage{graphicx}
\usepackage{bm}
\usepackage{tikz}
\usepackage{multirow}
\usepackage{color}

\graphicspath{ {./figs/} }

\title{Nonlocal multicontinua upscaling for multicontinua flow problems in fractured porous media}

\author{
Maria Vasilyeva \thanks{Institute for Scientific Computation, Texas A\&M University, College Station, TX 77843-3368 \& Department of Computational Technologies, North-Eastern Federal University, Yakutsk, Republic of Sakha (Yakutia), Russia, 677980. Email: {\tt vasilyevadotmdotv@gmail.com}.}
\and
Eric T. Chung \thanks{Department of Mathematics,
The Chinese University of Hong Kong (CUHK), Hong Kong SAR. Email: {\tt tschung@math.cuhk.edu.hk}.}
\and
Siu Wun Cheung \thanks{Department of Mathematic,
Texas A\&M University, College Station, Texas, USA}
\and
Yating Wang
\thanks{Department of Mathematics, Texas A\&M University, College Station, TX 77843-3368, USA.}
\and
Georgy Prokopev
\thanks{Multiscale model reduction laboratory, North-Eastern Federal University, Yakutsk, Republic of Sakha (Yakutia), Russia, 677980.}
}

\begin{document}

\maketitle

\begin{abstract}
Our goal of this paper is to develop a new upscaling method for multicontinua flow problems in fractured porous media. We consider a system of equations that describes flow phenomena with multiple flow variables defined on both matrix and fractures. To construct our upscaled model, we will apply the nonlocal multicontinua (NLMC) upscaling technique. The upscaled coefficients are obtained by using some multiscale basis functions, which are solutions of local problems defined on oversampled regions. For each continuum within a target coarse element, we will solve a local problem defined on an oversampling region obtained by extending the target element by few coarse grid layers, with a set of constraints which enforce the local solution to have mean value one on the chosen continuum and zero mean otherwise. The resulting multiscale basis functions have been shown to have good approximation properties. To illustrate the idea of our approach, we will consider a dual continua background model consisting of discrete fractures in two space dimensions, that is, we consider a system with three continua. We will present several numerical examples, and they show that our method is able to capture the interaction between matrix continua and discrete fractures on the coarse grid efficiently. 

\end{abstract}

\section{Introduction}

Multicontinuum models are common in many applications, for example, in geothermal reservoirs, oil and gas production, nuclear waste disposal etc.
In gas production from shale formation, we have highly heterogeneous and a complex mixture of organic matter, inorganic matter and multiscale fractures \cite{akkutlu2012multiscale, akkutlu2017multiscale}. The small scale heterogeneities in such unconventional reservoirs can be described by multicontinua models with additional lower dimensional discrete fracture for simulation of the flow in hydraulic fractures \cite{chung2017coupling, li2018multiscale}. Another example is the fractured vuggy reservoirs, where multicontinuum models are used for characterization of the complex interaction between vugges, fractures and porous matrix \cite{wu2011multiple, wu2007triple, yao2010discrete}.

Simulation of the flow problems in the fractured porous media is challenging due to multiple scales and high contrast of the properties.
The fracture network requires a special approach in the construction of a mathematical model and computational algorithms. Due to the high permeability of the fractures, they have a significant effect on the flow processes. Moreover,  fractures are different by the nature of their occurrence (naturally fractured media, faults, hydraulic fracturing technology).

Complex processes in the fractured porous media require some type of the model reduction \cite{houwu97, eh09, egh12, lunati2006multiscale, jenny2005adaptive}.
Typical approaches use an effective media property by the solution of the local problems in each coarse cell or representative volume \cite{sanchez1980non, talonov2016numerical}. Such approaches are not efficient for the case when each coarse grid blocks contains multiple important modes \cite{douglas1990dual, xu2001modeling, warren1963behavior, barenblatt1960basic}. For such problems the multicontinuum models are usually used. Specifically, in naturally fractured porous media, the fracture network is highly connected and the double porosity (dual continuum) approaches are employed. Such models are built for an ideal case and have a number of limitations. The interaction of the continuum in such models is determined by specifying the transfer functions between the matrix and the fractures. The definition of these transfer functions is a key issue for multicontinuum approach, since the construction of the exchange term is based on additional assumptions.
In the case of the large fractures, the fracture networks should be considered explicitly using Discrete Fracture Model (DFM) or Embedded Fracture Model (EFM) \cite{geothermal2018multiscale, ctene2016algebraic, karimi2003efficient, karimi2001numerical}. Accurate and explicit consideration of the fracture flow and complex interaction with porous matrix leads to the large system of equations and is computationally expensive.
Several multiscale methods for such problems are presented, for example, MsFEM \cite{hkj12, ctene2016algebraic, bosma2017multiscale, ctene2017projection} or GMsFEM \cite{akkutlu2015multiscale, chung2017coupling, akkutlu2018multiscale, li2018multiscale,CYH2016adaptive}.
Recently, the authors in \cite{chung2018nonfrac} proposed a new nonlocal multicontinua method (NLMC) for construction of the upscaled coarse grid model. In NLMC, one constructs multiscale basis functions which can capture complex matrix-fracture interaction by computing of the multiscale basis functions in a local oversampling domain. The construction of the multiscale space starts with the definition of the constraints that provide physical meanings for the coarse grid solution. Based on the constraints, one can obtain the required multiscale spaces by solving a constraint energy minimization problem in local domains. Moreover, since the local solutions are computed in an oversampled domain, the mass transfers between fractures and matrix become non-local, and the resulting upscaled model contains more effective properties of the flow problem. Recently, we extend this method for different applied problems \cite{chung2018nonfrac, chung2018nonfrac, vasilyeva2018nonlocal, vasilyeva2018constrained, vasilyeva2018upscaling}.

In this work, we consider the triple-continua model, where the background medium is considered to be dual continua with additional lower-dimensional fracture, and there are corresponding mass transfer between background medium and large fractures.
The mathematical model is described by a coupled mixed dimensional problem for simulation of the flow process in the fractured porous media with dual continuum background.
In order to reduce the size of the system and efficiently obtain solution of the presented problem, we propose a coarse grid approximation using nonlocal multicontinuum method and construct coupled multiscale basis functions.
The implementation is based on the open-source library FEniCS, where we use geometry objects and interface to the linear solvers \cite{logg2009efficient, logg2012automated}.

This paper is organized as follows.
In Section 2, we present problem formulation, where we consider a multicontinuum model with discrete fractures and concentrate on the triple continuum model. In Section 3, we consider fine grid approximation using finite volume method for the coupled system. Next, we construct upscaled model using NLMC method for coupled system of equation in Section 4. Finally, in Section 5, we construct an efficient numerical implementation and perform numerical investigation on two coarse grids and different number of the oversampling layers in local domain construction. We note that, using property of the constructed multiscale basis functions, we obtain a meaningful coarse grid upscaled model. We present results of the numerical simulations for two-dimensional model problem with dual continuum background medium and discrete fractures.

\section{Problem formulation}

\begin{figure}[h!]
  \centering
  \includegraphics[width=0.7 \textwidth]{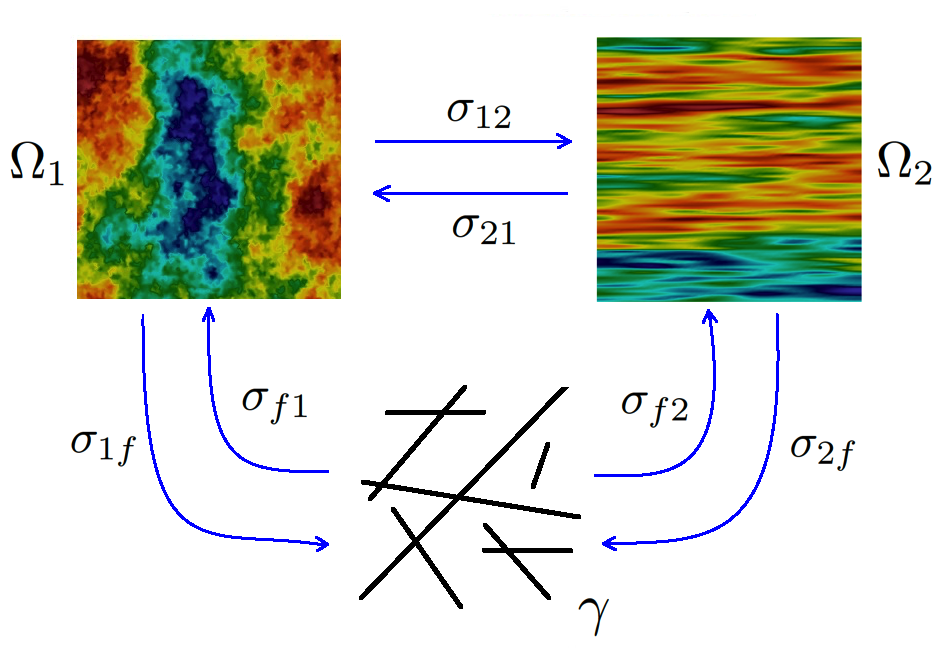}
  \caption{Illustration of the multicontinuum model on computational domain $\Omega_1 \cup \Omega_2 \cup \gamma$}
  \label{sch}
\end{figure}

Let $\Omega_1$ and $\Omega_2$ be the computational domains for first and second continua (Fig. \ref{sch}). For flow in large scale fractures, we consider a lower dimensional domain $\gamma = \cup_l \gamma_l$ \cite{martin2005modeling, formaggia2014reduced}.
The mass conservation and Darcy laws for the flow problem in  $\Omega_1 \cup \Omega_2 \cup \gamma$ indicates
\begin{equation}
\label{mm1}
\begin{split}
& c_{1} \frac{ \partial  p_1}{\partial t}
+ \nabla \cdot  q_{1}
+   \sigma_{1 2} (p_1 - p_2)
+   \sigma_{1 f} (p_1 - p_f) = f_1, \quad
x \in \Omega_1 \\
& q_{1} = - k_{1}  \nabla p_{1}, \quad  x \in \Omega_1, \\
& c_{1} \frac{ \partial p_2}{\partial t}
+ \nabla \cdot  q_{2}
+   \sigma_{2 1} (p_2 - p_1)
+   \sigma_{2 f} (p_2 - p_f)  = f_2, \quad
x \in \Omega_2 \\
& q_{2} = - k_{2}  \nabla p_{2}, \quad  x \in \Omega_2, \\
& c_{f} \frac{ \partial  p_f}{\partial t}
+ \nabla_\gamma \cdot q_{f}
+   \sigma_{f 1} (p_f - p_1)
+   \sigma_{f 2} (p_f - p_2)  = f_f, \quad
x \in \gamma \\
& q_{f} = - k_{f}  \nabla_\gamma p_{f}, \quad x \in \gamma,
\end{split}
\end{equation}
where $q_{\alpha}$, $p_{\alpha}$, $f_{\alpha}$ are the flux, pressure and source term for the $\alpha$ continuum with $\alpha = 1, 2, f$, where  subindices $1, 2$ are related for first and second continuums, and subindex $f$ is related to the lower-dimensional fracture model.
Here $k_{\alpha} = \kappa_{\alpha}/\mu$, $\mu$ is the fluid viscosity, $\kappa_{\alpha}$ and $c_{\alpha}$  are the permeability and compressibility of the $\alpha$ continuum ($\alpha = 1, 2, f$).
System of equations \eqref{mm1} are coupled by the mass exchange terms between continua $\sigma_{\alpha \beta}$ ($\sigma_{\alpha \beta} = \sigma_{\beta \alpha}$), where $\alpha, \beta = 1, 2$ and $\beta \neq \alpha$.

 Substituting the Darcy's equatio  into the mass conservation equation, we obtain following system of equations for $p_1$, $p_2$ and $p_f$
\begin{equation}
\label{mm2}
\begin{split}
& c_{1} \frac{ \partial  p_1}{\partial t}
- \nabla \cdot ( k_{1}  \nabla p_{1})
+   \sigma_{1 2} (p_1 - p_2)
+   \sigma_{1 f} (p_1 - p_f) = f_1, \quad
x \in \Omega_1 \\
& c_{1} \frac{ \partial p_2}{\partial t}
- \nabla \cdot (k_{2}  \nabla p_{2})
+   \sigma_{2 1} (p_2 - p_1)
+   \sigma_{2 f} (p_2 - p_f)  = f_2, \quad
x \in \Omega_2 \\
& c_{f} \frac{ \partial  p_f}{\partial t}
- \nabla_\gamma \cdot (k_{f}  \nabla_\gamma p_{f})
+   \sigma_{f 1} (p_f - p_1)
+   \sigma_{f 2} (p_f - p_2)  = f_f, \quad
x \in \gamma.
\end{split}
\end{equation}
We consider system of equations \eqref{mm2} with the homogeneous Neumann boundary conditions and some given initial conditions, $p_1 = p_2 = p_f = p_0$ for $t = 0$.

Note that, we can write the flow problem using a general multicontinuum model 
\begin{equation}
\label{mm3}
\begin{split}
& c_{\alpha} \frac{ \partial  p_{\alpha}}{\partial t}
+ \nabla \cdot  q_{\alpha}
+  \sum_{\beta \neq \alpha} \sigma_{\alpha \beta} (p_{\alpha} - p_{\beta})  = f_{\alpha}, \\
& q_{\alpha} = - k_{\alpha}  \nabla p_{\alpha},
\end{split}
\end{equation}
or we can write
\begin{equation}
\label{mm4}
c_{\alpha} \frac{ \partial  p_{\alpha}}{\partial t}
- \nabla \cdot  ( k_{\alpha}  \nabla p_{\alpha} )
+  \sum_{\beta \neq \alpha} \sigma_{\alpha \beta} (p_{\alpha} - p_{\beta})  = f_{\alpha},
\end{equation}
where $\alpha, \beta = 1,2,...,N$ and $N$ is the number of continua.

\section{Approximation on the fine grid}

Let $\mathcal{T}^h$ denote a triangulation of the domain $\Omega = \Omega_1 = \Omega_2$
and $\cup_j \gamma_j$ represent fractures, where  $j = \overline{1, N_{frac}}$ and $N_{frac}$ is the number of discrete fractures.

For approximation, we used finite volume method with discrete or embedded fracture models on the fine grid for ($p_1$, $p_2$, $p_f$)
\begin{equation}
\label{app1}
\begin{split}
& c_1 \frac{ p_{1, i} - \check{p}_{1, i} }{\tau} |\varsigma_i|
 + \sum_{j}  T_{1,ij}  (p_{1, i} - p_{1, j})
 +  \sigma_{12,ii} (p_{1, i} - p_{2, i} )
 +  \sigma_{1f,in} (p_{1, i} - p_{f, n} )
 =  f_1   |\varsigma_i|, \quad \forall i = 1, N^1_f \\
 & c_2 \frac{ p_{2, i} - \check{p}_{2, i} }{\tau} |\varsigma_i|
 + \sum_{j}  T_{2,ij}  (p_{2, i} - p_{2, j})
 +  \sigma_{12,ii} (p_{2, i} - p_{1, i} )
 +  \sigma_{2f,in} (p_{2, i} - p_{f, n} )
 =  f_2   |\varsigma_i|, \quad \forall i = 1, N^2_f \\
& c_f \frac{ p_{f, l} - \check{p}_{f, l}}{\tau}  |\iota_l|
+ \sum_{n} T_{f,ln} (p_{f, l} - p_{f, n})
+ \sigma_{il} (p_{f, n} - p_{1, i} )
+ \sigma_{in} (p_{f, n} - p_{2, i} )
 =  f_f  |\iota_l|, \quad \forall l = 1, N^f_f
\end{split}
\end{equation}
where
$T_{1,ij} = k_1 |E_{ij}|/d_{ij}$, $T_{2,ij} = k_2 |E_{ij}|/d_{ij}$ ($|E_{ij}|$ is the length of facet between cells $\varsigma_i$ and $\varsigma_j$, $d_{ij}$ is the distance between midpoint of cells $\varsigma_i$ and $\varsigma_j$),
$T_{f,ln} = b_f/d_{ln}$ ($d_{ln}$ is the distance between points $l$ and $n$),
$N^m_f$ is the number of cells in $\mathcal{T}_h$,
$N^f_f$ is the number of cell related to fracture mesh $\mathcal{E}_{\gamma}$,
$\sigma_{il} = \sigma$ if $\mathcal{E}_{\gamma} \cap \partial \varsigma_i = \iota_l$ and zero else.
We use implicit scheme for time discretization and $(\check{p}_1, \check{p}_2,  \check{p}_f)$ is the solution from previous time step and $\tau$ is the given time step.

We can write similar approximation for multicontinuum model
\begin{equation}
\label{app2}
 c_{\alpha} \frac{ p_{{\alpha}, i} - \check{p}_{{\alpha}, i} }{\tau} |\varsigma_i|
 + \sum_{j}  T_{{\alpha},ij}  (p_{{\alpha}, i} - p_{{\alpha}, j})
 + \sum_{{\alpha \neq \beta}} \sigma_{\alpha \beta,ii} (p_{{\alpha}, i} - p_{\beta, i} )
 =  f_{\alpha}   |\varsigma_i|, \quad \forall i = 1, N^{\alpha}_f,
\end{equation}
where $\alpha = 1,2...$

In the matrix form, we have following system for $p = (p_1, p_2, p_f)$
\begin{equation}
\label{app3}
M \frac{p - \check{p}}{\tau} + (A + Q) p = F,
\end{equation}
where
\begin{equation}
\label{app4}
M = \begin{pmatrix}
M_1 & 0 & 0 \\
0 & M_2 & 0 \\
0 & 0 & M_f
\end{pmatrix},  \quad
A = \begin{pmatrix}
A_1 & 0 & 0 \\
0 & A_2 & 0 \\
0 & 0 & A_f
\end{pmatrix},
\end{equation}
\begin{equation}
\label{eq:fa-4}
Q = \begin{pmatrix}
Q_{12} + Q_{1f} & -Q_{12} & -Q_{1f} \\
-Q_{12} & Q_{12} + Q_{2f} & -Q_{2f} \\
-Q_{1f}   & -Q_{2f} & Q_{1f} + Q_{2f}
\end{pmatrix},
\end{equation}
and
\[
M_{1} = \{m^{1}_{ij}\}, \quad
m^{1}_{ij} =
\left\{\begin{matrix}
c_{1} |\varsigma_i|  & i = j, \\
0 & i \neq j
\end{matrix}\right. , \quad
M_{2} = \{m^{2}_{ij}\}, \quad
m^{2}_{ij} =
\left\{\begin{matrix}
c_{2} |\varsigma_i|  & i = j, \\
0 & i \neq j
\end{matrix}\right. ,
\]\[
M_f = \{m^f_{ln}\}, \quad
m^f_{ln} =
\left\{\begin{matrix}
 c_f |\iota_l| / \tau & l = n, \\
0 & l \neq n
\end{matrix}\right. ,\quad
Q_{\alpha \beta} = \{q_{\alpha \beta, il}\}, \quad
q_{\alpha \beta,il} =
\left\{\begin{matrix}
\sigma_{\alpha \beta} & i = l, \\
0 & i \neq l
\end{matrix}\right. ,
\]
where
$A_{1} = \{T_{1, ij}\}$, $A_{2} = \{T_{2, ij}\}$,   $A_f = \{W_{ln}\}$,
$F_{1} = \{f^{1}\}$, $F_{2} = \{f^{2}\}$, $f^m_i = f_m |\varsigma_i|$,
$F_f = \{f^f_l\}$, $f^f_l = f_f |\iota_l|$ and size of fine-grid system is $N_f = N^1_f + N^2_f + N^f_f$.

\section{Non-local multi-continua upscaling on the coarse grid}

Next, we  construct an accurate approximation of the coupled system of equations using the NLMC approach.
To construct multiscale basis functions, one can solve local problems in local domain satisfying the flow equations and subject to some constraints. We note the meaning of the constraints is that the local solution has zero mean in other continua except the one for which it is formulated for.
The resulting multiscale basis functions have spatial decay property in local domains and can separate continua.
The resulting basis functions will be used in the construction of the upscaled model.

\begin{figure}[h!]
\centering
\includegraphics[width=0.6\textwidth]{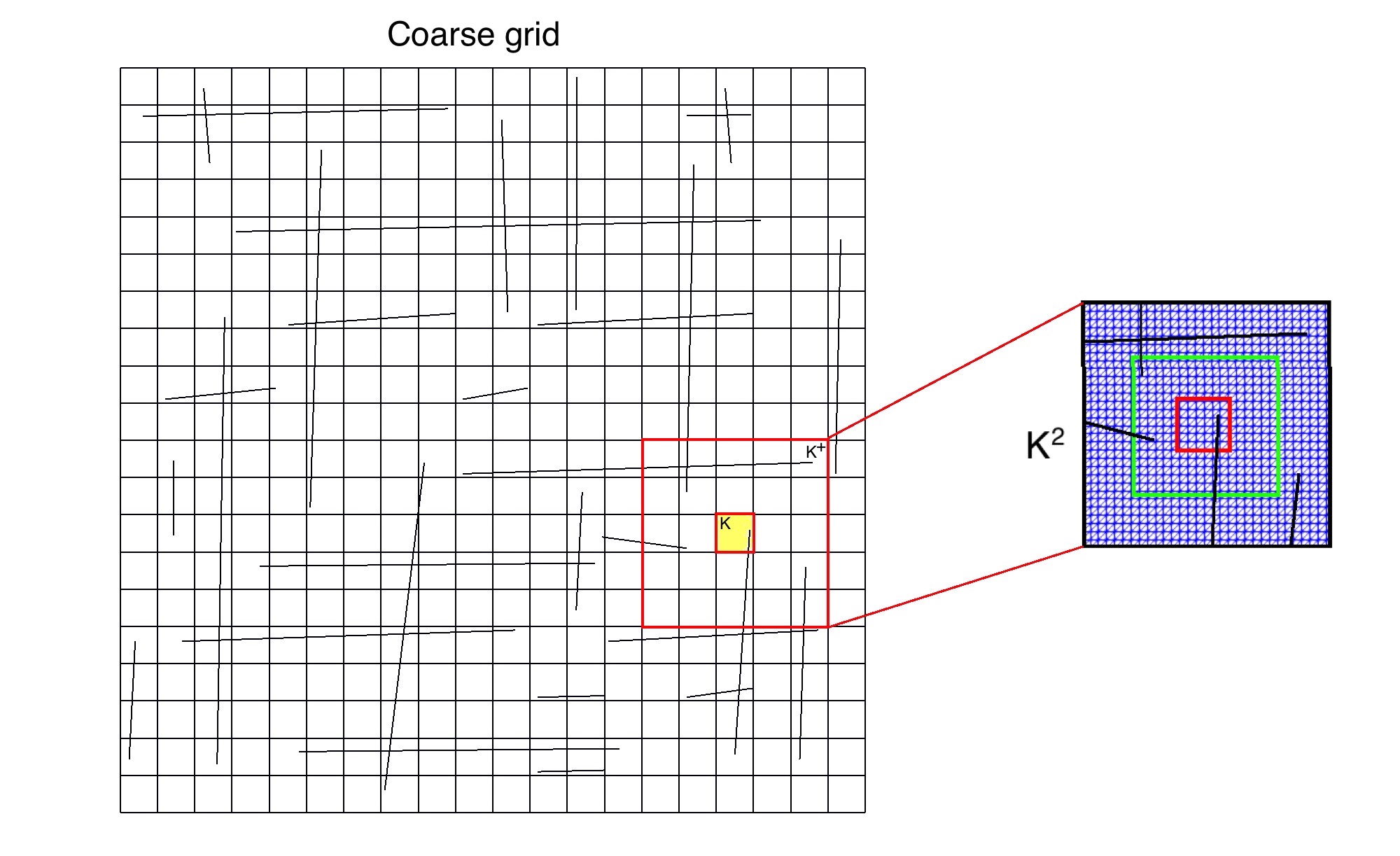}
\caption{Illustration of the coarse mesh, fracture distribution and local domain. Coarse mesh $20 \times 20$.  Local domain with two-oversampling layers,  $K^2$}
\label{fig:sch1}
\end{figure}

Let $K^+_i$ be an oversampled region for the coarse cell $K_i$ (see Figure \ref{fig:sch1}) obtained by enlarging $K_i$ by several coarse cell layers.
We will construct a set of basis functions, whose supports are $K_i^+$.
Each of these basis functions is related to the matrix component in $K_i$ as well as each fracture network within $K_i$.
For fractures, we denote $\gamma = \cup_{l = 1}^L \gamma^{(l)}$, where $\gamma^{(l)}$ denotes the $l$-th fracture network and $L$ is the total number of fracture networks.
We also write $\gamma^{(l)}_j = K_j \cap \gamma^{(l)}$ as the fracture inside coarse cell $K_j$, and $L_j$ is the number of fracture networks in $K_j$.
For each $K_i$, we will therefore obtain $L_j+2$ basis functions: two for $K_i$ (for dual continuum background model) and one for each $\gamma^{(l)}_i$.

For the $i$-th coarse element $K_i$ and the $l$-th continuum within $K_i$, we will obtain $3$ basis functions $\psi^{i,l} = (\psi^{i,l}_1,\psi^{i,l}_2,\psi^{i,l}_f)$, whose supports are $K_i^+$. 
We denote $a(\psi,v)$ as the variational form of the spatial differential operator restricted on $K_i^+$. Then we can will find the basis function $\psi^{i,l} = (\psi^{i,l}_1,\psi^{i,l}_2,\psi^{i,l}_f)$
such that
\begin{equation*}
a(\psi^{i,l},v) = \langle \mu,v\rangle
\end{equation*}
for all suitable test functions $v$, and subject to the following constraints
\begin{equation*}
\int_{\gamma_j^{(m)}} \psi^{i,l}_{\alpha} = 0, \quad \forall (j,m) \ne (i,l), \text{ and } \alpha = 1,2,f
\end{equation*}
and
\begin{equation*}
\int_{\gamma_i^{(l)}} \psi^{i,l}_{\alpha} = 1, \;\text{ and for one fixed } \alpha \text{ and zero otherwise}.
\end{equation*}
We note that the above is a set of constraints so that the resulting function has mean value one on the coarse cell $K_i$ for current continuum, and has mean value zero on all other coarse cells within $K_i^+$ and on all coarse cells for another continua.

Below, we will write the above formulation in a more concrete setting. In particular, we have 
\begin{equation}
\label{eq:basis}
\begin{pmatrix}
A_1+ Q_{12}+ Q_{1f} & -Q_{12} & -Q_{1f} & B^T_1 & 0 & 0 \\
-Q_{12} & A_2+ Q_{12}+ Q_{2f} & -Q_{2f} & 0 & B^T_2 & 0 \\
-Q_{1f} & -Q_{2f} & A_f + Q_{1f} + Q_{2f} & 0  & 0 &  B^T_f \\
B_1 & 0 & 0 & 0  & 0 & 0 \\
0 & B_2 & 0 & 0  & 0  & 0\\
0 & 0 & B_f & 0  & 0  & 0\\
\end{pmatrix}
\begin{pmatrix}
\psi_1 \\
\psi_2 \\
\psi_f \\
\mu_1 \\
\mu_2 \\
\mu_f \\
\end{pmatrix} =
\begin{pmatrix}
0 \\
0 \\
0 \\
F_1 \\
F_2 \\
F_f \\
\end{pmatrix}
\end{equation}
with  zero Dirichlet boundary conditions on $\partial K^+_i$ for $\psi_{\alpha}$, $\alpha = 1,2,f$.
We remark that $(\psi_1, \psi_2, \psi_f)$ denotes each of the basis functions that satisfy the above constraints.
Note that we used Lagrange multipliers $\mu_{\alpha}$ to impose the constraints in the multiscale basis construction,
and that the matrices $B_{\alpha}$ are the mean value operators. 

To construct multiscale basis function with respect to each continuum,  we set $F_{\alpha} = \delta_{i,j}$ and $F_f = 0$ for $\alpha = 1,2$. For multiscale basis function with respect to the $l$-th fracture network, we set  $F_1 = F_2 =  0$ and $F_f = \delta_{i,j}\delta_{m,l}$. In Figure \ref{fig:sch2}, we depict a multiscale basis functions for each continuum in oversampled region $K^+_i = K^2_i$ (two oversampling coarse cell layers) in a $20 \times 20$ coarse mesh.
Combining these multiscale basis functions, we obtain the following multiscale space
\[
V_{ms} = \text{span} \{ \psi^{i,l}_{\alpha}, \quad \alpha = 1,2,f, \quad i = \overline{1,N_c}, \quad l = \overline{1, L_i+2} \}
\]
and the projection matrix
\[
R = \begin{pmatrix}
R_{11} & R_{12}& R_{1f} \\
R_{12} & R_{22}& R_{2f} \\
R_{1f} & R_{2f}  & R_{ff}
\end{pmatrix},
\]\[
R_{\alpha \alpha}^T = \left[
\psi^{0,\alpha}_{{\alpha}}, \psi^{1,\alpha}_{{\alpha}} \ldots \psi^{N_c,\alpha}_{\alpha}
 \right], \quad
R_{\alpha \beta}^T = \left[
\psi^{0,\alpha+1}_{\alpha} \ldots \psi^{0,L_{\alpha+1}}_{\alpha},
\psi^{1,\alpha+1}_{\alpha} \ldots \psi^{1,L_{\alpha+2}}_{\alpha},
\ldots,
\psi^{N_c,\alpha+1}_{\alpha}  \ldots \psi^{N_c,L_{\alpha+N_c}}_{\alpha}
 \right],
\]\[
R_{mf}^T = \left[
\psi^{0,0}_f, \psi^{1,0}_f \ldots \psi^{N_c,0}_f
 \right], \quad
R_{fm}^T = \left[
\psi^{0,1}_m \ldots \psi^{0,L_0}_m,
\psi^{1,1}_m \ldots \psi^{1,L_1}_m,
\ldots,
\psi^{N_c,1}_m  \ldots \psi^{N_c,L_{N_c}}_m
 \right],
 \]
where  $\alpha,\beta = 1,2$.

Therefore, the resulting upscaled coarse grid model reads
\begin{equation}
\label{t-nlmc2}
\bar{M} \frac{\bar{p}^{n+1} - \bar{p}^n}{\tau} + \bar{A} \bar{p}^{n+1} = \bar{F},
\end{equation}
where
$\bar{A} = R A R^T$, $\bar{F} = R F$ and
$\bar{p} = (\bar{p}_1, \bar{p}_2, \bar{p}_f)$.
We remark that $\bar{p}_1$, $\bar{p}_2$  and $\bar{p}_f$ are the average cell solution on coarse grid element for background matrices and for fracture media.
That is, each component of $\bar{p}_m$ corresponds to the mean value of the solution on each coarse cell.
Moreover, each component of $\bar{p}_f$ corresponds to the mean value of the solution on each fracture network with a coarse cell.

As an approximation, we can use diagonal mass matrix directly calculated on the coarse grid
\[
\bar{M} =
\begin{pmatrix}
\bar{M}_] & 0 & 0 \\
0 & \bar{M}_2 & 0 \\
0 & 0 & \bar{M}_f \\
\end{pmatrix}, \quad
\bar{Q} =
\begin{pmatrix}
\bar{Q}_{12} + \bar{Q}_{1f}  &  - \bar{Q}_{12} &  - \bar{Q}_{1f} \\
-\bar{Q}_{12} & \bar{Q}_{12} + \bar{Q}_{2f} &  - \bar{Q}_{2f} \\
-\bar{Q}_{1f} &  - \bar{Q}_{2f} & \bar{Q}_{1g} + \bar{Q}_{2f}  \\
\end{pmatrix}, \quad
\bar{F} =
\begin{pmatrix}
\bar{F}_1 \\
\bar{F}_2 \\
\bar{F}_f \\
\end{pmatrix},
\]
where
$\bar{M}_{\alpha} = \text{diag}\{ a_{\alpha} |K_i| \}$ (${\alpha} = 1,2$),
$\bar{M}_f = \text{diag}\{ a_f |\gamma_i| \}$,
$\bar{Q}_{\alpha \beta} = \text{diag}\{ \sigma_{\alpha \beta} |\gamma_i| \}$ ($\alpha, \beta = 1,2,f$)
and for the right-hand side vector
$\bar{F}_{\alpha} = \{ f_{\alpha} |K_i| \}$ (${\alpha} = 1,2$),
$\bar{F}_f = \{ f_f |\gamma_i| \}$.
We remark that the matrix $A$ is non-local and provide good approximation due to the coupled multiscale basis construction.

\section{Numerical results}

\begin{figure}[h!]
\centering
\includegraphics[width=0.35 \textwidth]{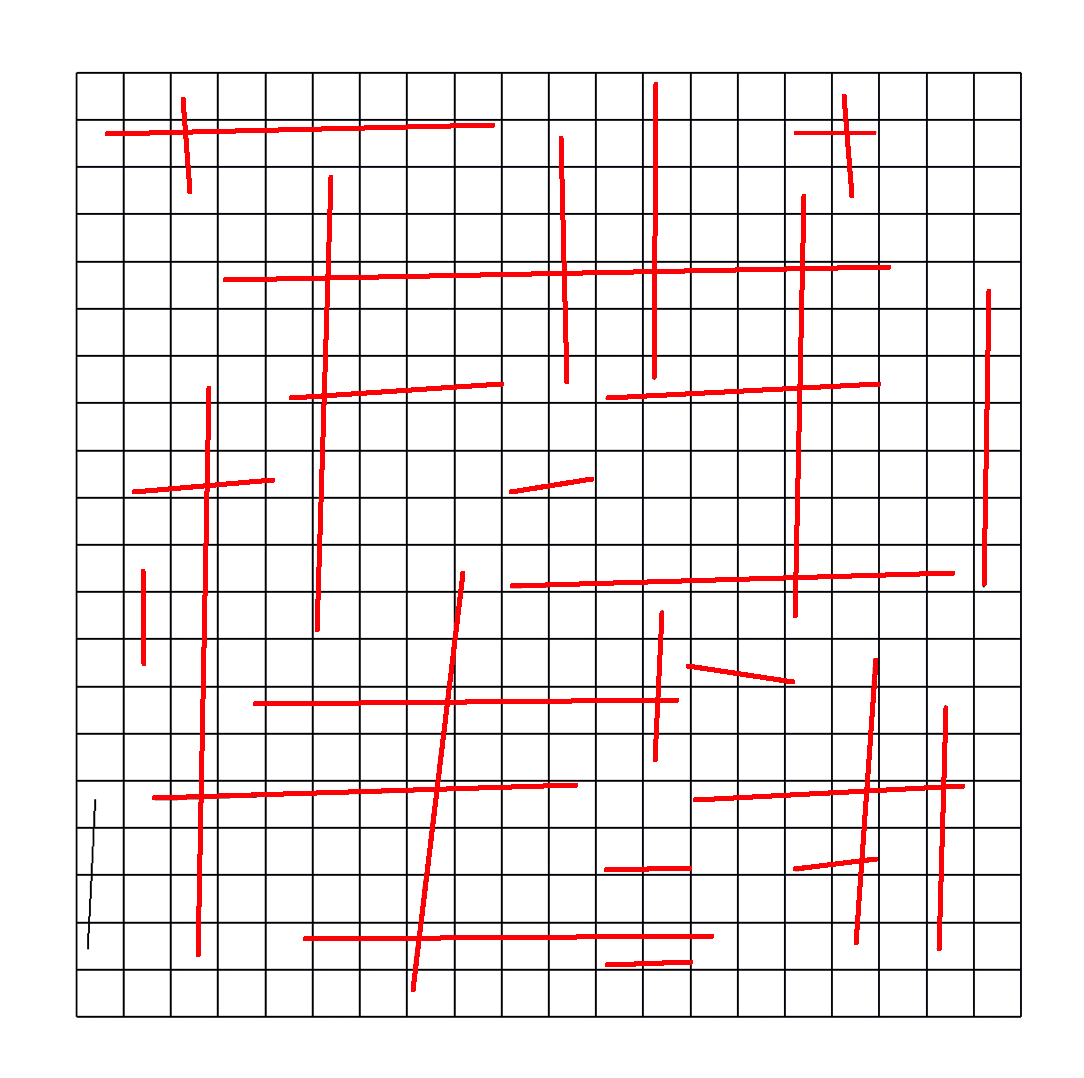}
\includegraphics[width=0.35 \textwidth]{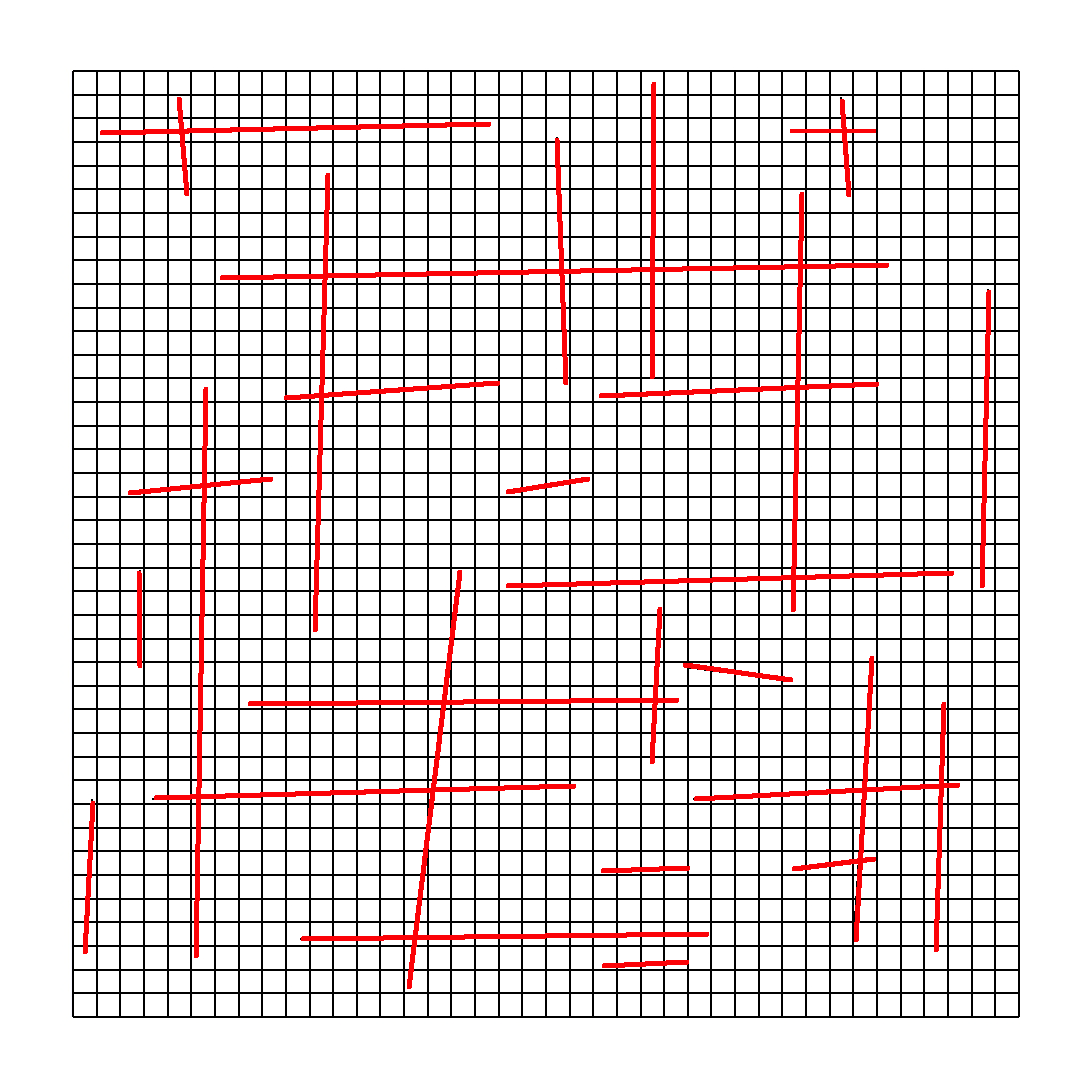}
\caption{Coarse mesh $20 \times 20$ with 400 cells (left) and $40 \times 40$ with 1600 cells (right). }
\label{fig:mesh}
\end{figure}

We consider a two-dimensional model problem with dual continuum background medium and discrete fractures. We present a numerical result for the proposed multiscale method in computational domain $\Omega = [0, 1] \times [0, 1]$.
We construct an efficient numerical implementation and perform numerical investigation on two coarse grids ($20 \times 20$ and $40 \times 40$) and different number of the oversampling layers in local domain construction, $K^s$ with $s = 1,2,3,4$. The coarse grid is uniform.
In Figure \ref{fig:mesh}, we depict the $20 \times 20$ and $40 \times 40$ coarse grids, where the fracture lines are depicted with red color.
The fractures are approximated on the fine grid using embedded fracture model (EFM) and fine grid is $120 \times 120$ (14400 quadrilateral cells). For fracture fine grid, we have $1042$ cells. Therefore, the fine grid system has 29842 degrees of freedom, $DOF_f = 29 842$.

\begin{figure}[h!]
\centering
\includegraphics[width=0.7\textwidth]{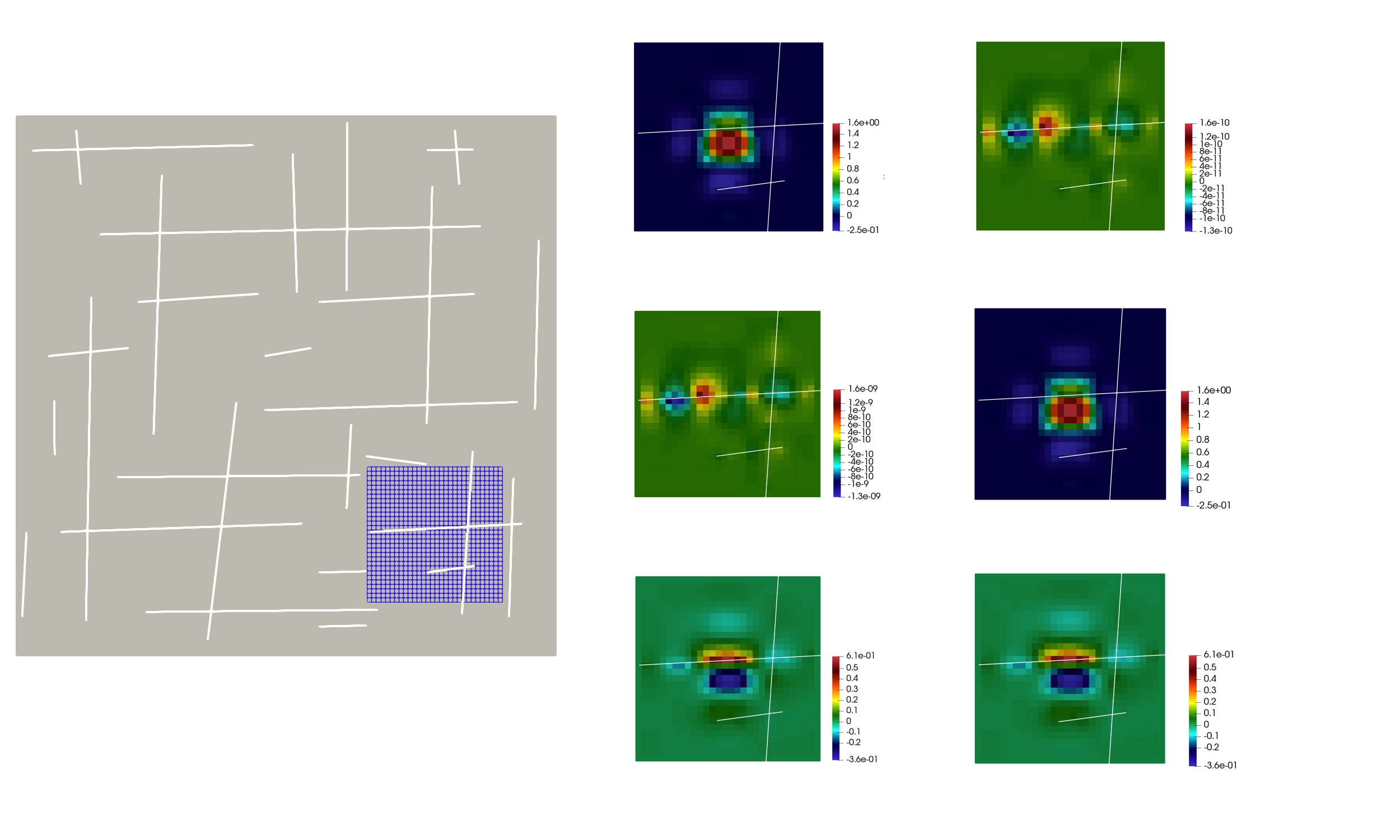}
\caption{Multiscale basis functions on mesh $20 \times 20$ for local domain $K^2$ for triple-continuum model (dual continuum background and fracture lines). \textit{Test 1}}
\label{fig:sch2}
\end{figure}

\begin{figure}[h!]
\centering
\includegraphics[width=0.7\textwidth]{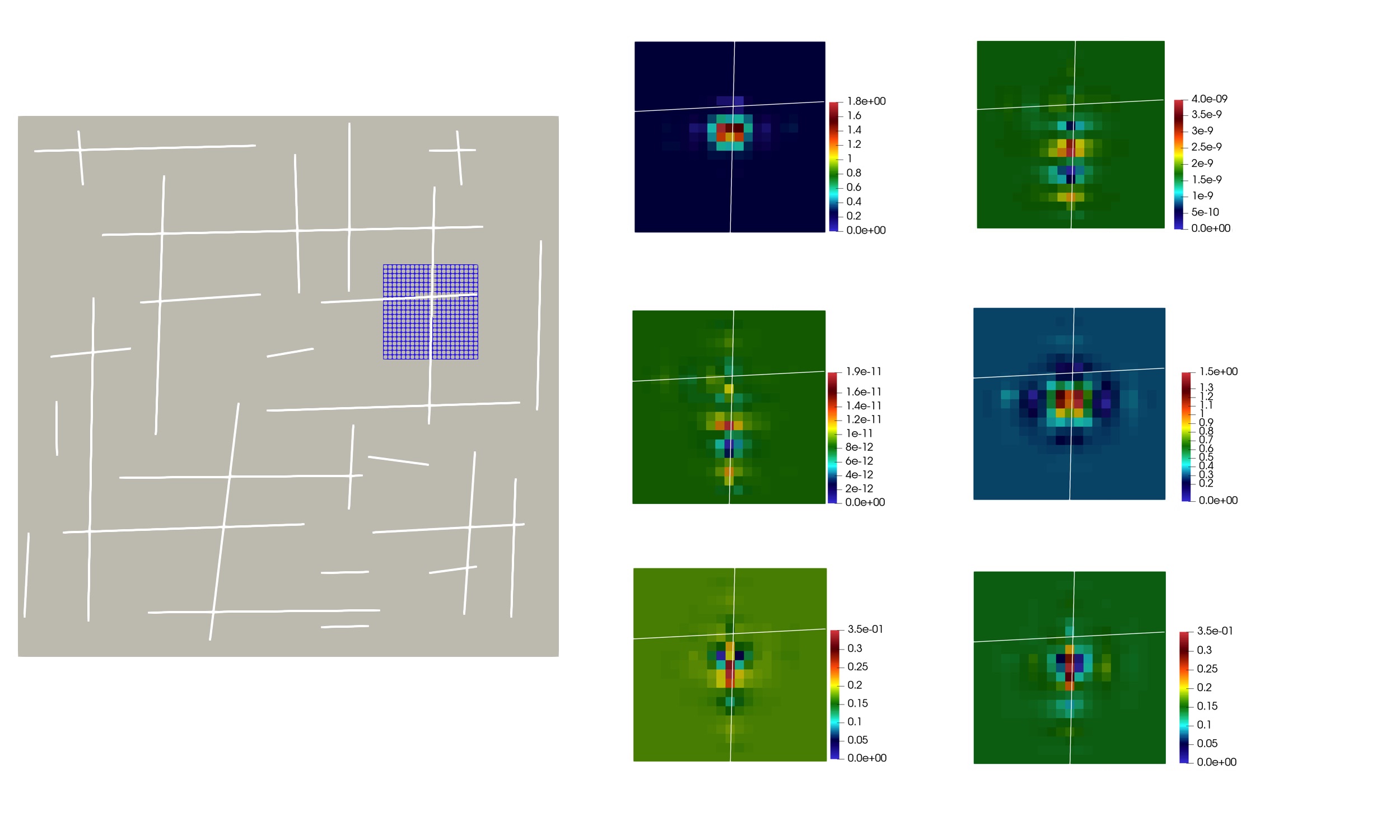}
\caption{Multiscale basis functions on mesh $40 \times 40$ for local domain $K^2$ for triple-continuum model (dual continuum background and fracture lines). \textit{Test 2}}
\label{fig:sch3}
\end{figure}

\begin{figure}[h!]
\centering
\includegraphics[width=0.32 \textwidth]{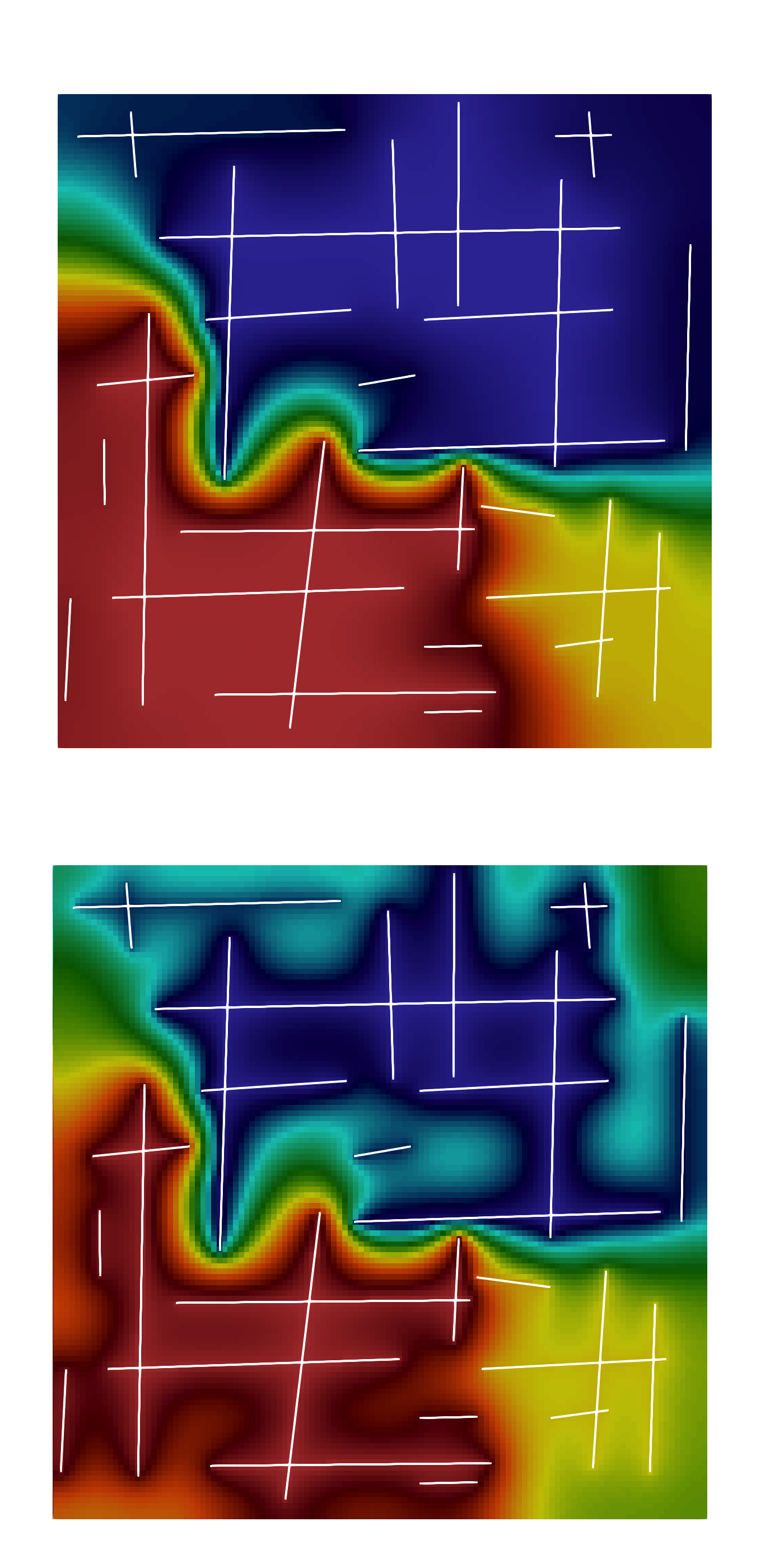}
\includegraphics[width=0.32 \textwidth]{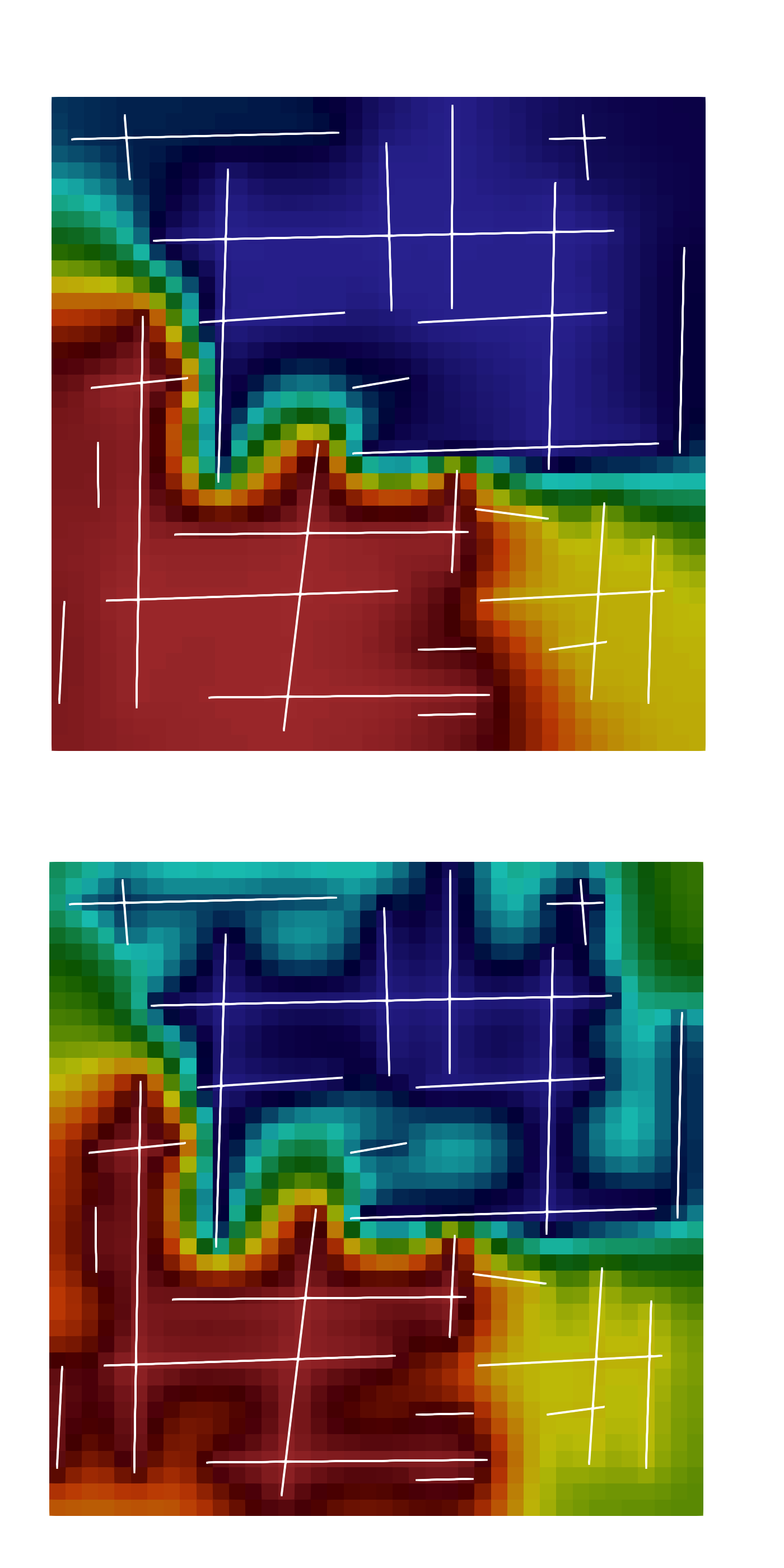}
\includegraphics[width=0.32 \textwidth]{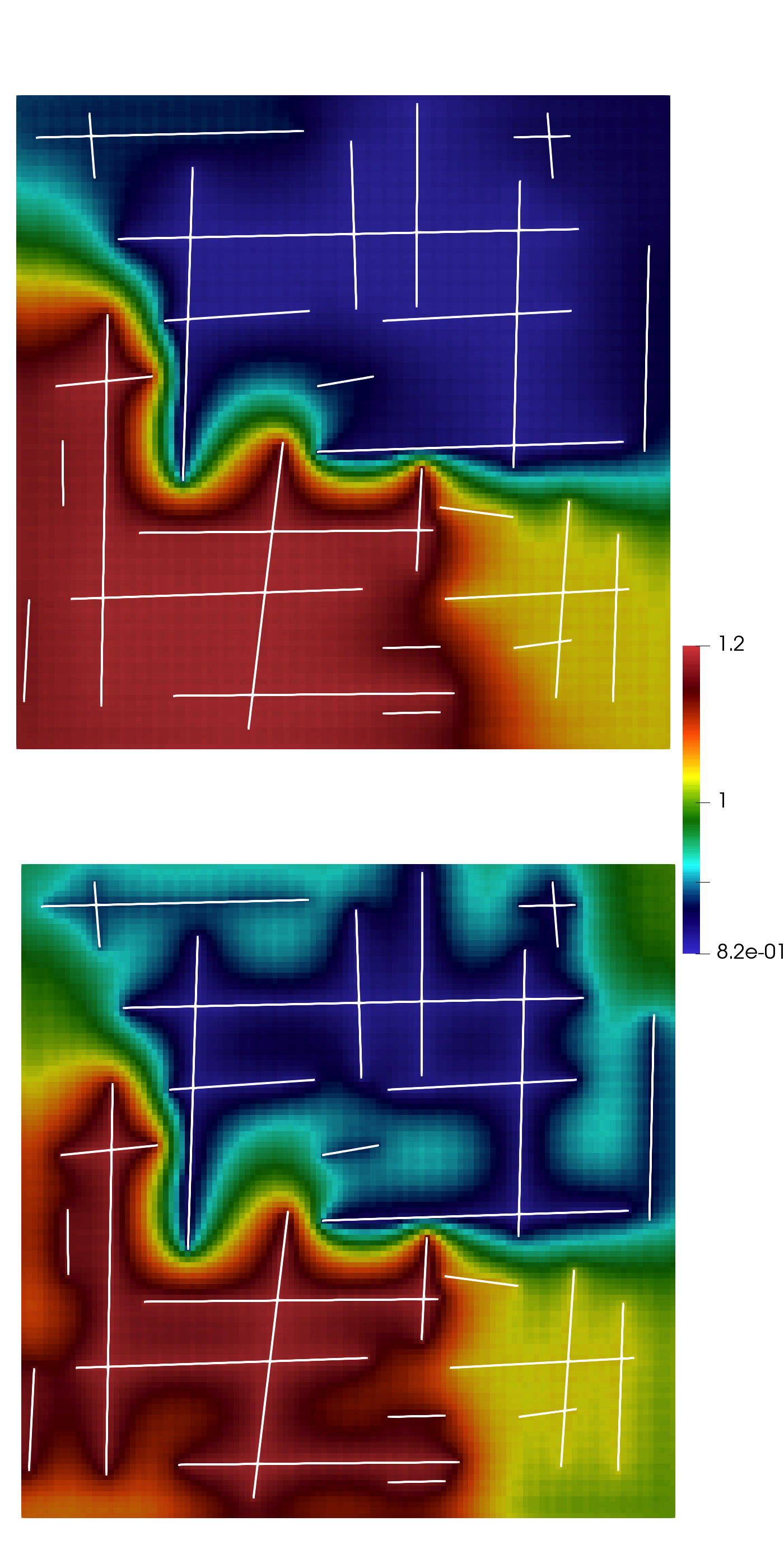}
\caption{Numerical results for pressure at final time for \textit{Test 1}.
Top: first continuum. Bottom: second continuum.
First column: fine-scale solution, $DOF_f = 29 842$.
Second column: coarse grid solution using NLMC, $DOF_c = 3 565$.
Third column: reconstructed fine grid solution using NLMC.
}
\label{fig:u1}
\end{figure}

\begin{table}[h!]
\centering
\begin{tabular}{ |c | c | c | c | c | }
\hline
$K^s$ &
$s = 1$&
$s = 2$&
$s = 3$&
$s = 4$ \\
\hline
\multicolumn{5}{|c|}{Coarse grid: $20 \times 20$} \\
\hline
$m = 5$ 		& 0.521 & 0.126 &	0.123 &	 0.124  \\ \hline
$m = 15$ 	& 0.520 & 0.133 &	0.132 &	 0.131 \\ \hline
$m = 25$ 	& 1.072 & 0.112 &	0.103 &	 0.102 \\ \hline
$m = 35$ 	& 1.678 & 0.114 &	0.076 &	 0.073 \\ \hline
$m = 50$ 	& 2.533 & 0.154 &	0.061 &	 0.054 \\ \hline
\multicolumn{5}{|c|}{Reconstructed fine grid} \\
\hline
$m = 5$ 		& 15.216 &	3.486 & 0.785 &	0.265  \\ \hline
$m = 15$ 	& 15.224 &	3.481 & 0.777 &	0.232  \\ \hline
$m = 25$ 	& 15.255 &	3.474 & 0.773 &	0.206  \\ \hline
$m = 35$ 	& 15.315 &	3.469 & 0.770 &	0.189  \\ \hline
$m = 50$ 	& 15.443 &	3.467 & 0.769 &	0.181  \\ \hline
\end{tabular}
\,\,\,
\begin{tabular}{ |c | c | c | c | c | }
\hline
$K^s$ &
$s = 1$&
$s = 2$&
$s = 3$&
$s = 4$ \\
\hline
\multicolumn{5}{|c|}{Coarse grid: $20 \times 20$} \\
\hline
$m = 5$ 		& 0.303 & 0.069 &	0.057 &	 0.057 \\ \hline
$m = 15$ 	& 0.850 & 0.109 &	0.058 &	 0.055 \\ \hline
$m = 25$ 	& 1.324 & 0.135 &	0.044 &	 0.039 \\ \hline
$m = 35$ 	& 1.712 & 0.161 &	0.034 &	 0.025 \\ \hline
$m = 50$ 	& 2.132 & 0.196 &	0.026 &	 0.014 \\ \hline
\multicolumn{5}{|c|}{Reconstructed fine grid} \\
\hline
$m = 5$ 		& 15.204 &	3.470 & 0.769 &	0.181  \\ \hline
$m = 15$ 	& 15.230 &	3.466 & 0.768 &	0.180  \\ \hline
$m = 25$ 	& 15.323 &	3.465 & 0.767 &	0.176  \\ \hline
$m = 35$ 	& 15.457 &	3.644 & 0.766 &	0.173  \\ \hline
$m = 50$ 	& 15.679 &	3.464 & 0.765 &	0.172  \\ \hline
\end{tabular}
\caption{Relative errors on the coarse grid and reconstructed fine-grid solutions for \textit{Test 1}. Coarse mesh: $20 \times 20$. $DOF_c = 993$ and $DOF_f = 29842$.
Left: first continuum. Right: second continuum. }
\label{err1a}
\end{table}

\begin{table}[h!]
\centering
\begin{tabular}{ |c | c | c | c | c | }
\hline
$K^s$ &
$s = 1$&
$s = 2$&
$s = 3$&
$s = 4$ \\
\hline
\multicolumn{5}{|c|}{Coarse grid: $40 \times 40$} \\
\hline
$m = 5$ 		& 0.192 & 0.046 &	0.039 &	 0.039  \\ \hline
$m = 15$ 	& 0.531 & 0.097 &	0.038 &	 0.038 \\ \hline
$m = 25$ 	& 0.836 & 0.170 &	0.029 &	 0.028 \\ \hline
$m = 35$ 	& 1.124 & 0.243 &	0.022 &	 0.019 \\ \hline
$m = 50$ 	& 1.496 & 0.336 &	0.023 &	 0.012 \\ \hline
\multicolumn{5}{|c|}{Reconstructed fine grid} \\
\hline
$m = 5$ 		& 13.919 &	3.098 & 0.657 &	0.150  \\ \hline
$m = 15$ 	& 13.941 &	3.105 & 0.657 &	0.148  \\ \hline
$m = 25$ 	& 13.961 &	3.111 & 0.657 &	0.146  \\ \hline
$m = 35$ 	& 13.982 &	3.117 & 0.657 &	0.144  \\ \hline
$m = 50$ 	& 14.013 &	3.126 & 0.657 &	0.143  \\ \hline
\end{tabular}
\,\,\,
\begin{tabular}{ |c | c | c | c | c | }
\hline
$K^s$ &
$s = 1$&
$s = 2$&
$s = 3$&
$s = 4$ \\
\hline
\multicolumn{5}{|c|}{Coarse grid: $40 \times 40$} \\
\hline
$m = 5$ 		& 0.303 & 0.046 &	0.013 &	 0.016 \\ \hline
$m = 15$ 	& 0.850 & 0.166 &	0.008 &	 0.014 \\ \hline
$m = 25$ 	& 1.324 & 0.279 &	0.016 &	 0.008 \\ \hline
$m = 35$ 	& 1.712 & 0.376 &	0.026 &	 0.003 \\ \hline
$m = 50$ 	& 2.132 & 0.483 &	 0.038 &	 0.001 \\ \hline
\multicolumn{5}{|c|}{Reconstructed fine grid} \\
\hline
$m = 5$ 		& 13.924 &	3.100 & 0.655 &	0.143  \\ \hline
$m = 15$ 	& 13.957 &	3.108 & 0.656 &	0.143  \\ \hline
$m = 25$ 	& 13.995 &	3.118 & 0.656 &	0.142  \\ \hline
$m = 35$ 	& 14.035 &	3.129 & 0.656 &	0.142  \\ \hline
$m = 50$ 	& 14.087 &	3.144 & 0.657 &	0.142  \\ \hline
\end{tabular}
\caption{Relative errors on the coarse grid and reconstructed fine-grid solutions for \textit{Test 1}. Coarse mesh: $40 \times 40$.  $DOF_c  = 3565$ and $DOF_f = 29842$.
Left: first continuum. Right: second continuum. }
\label{err1}
\end{table}

\begin{figure}[h!]
\centering
\includegraphics[width=0.32 \textwidth]{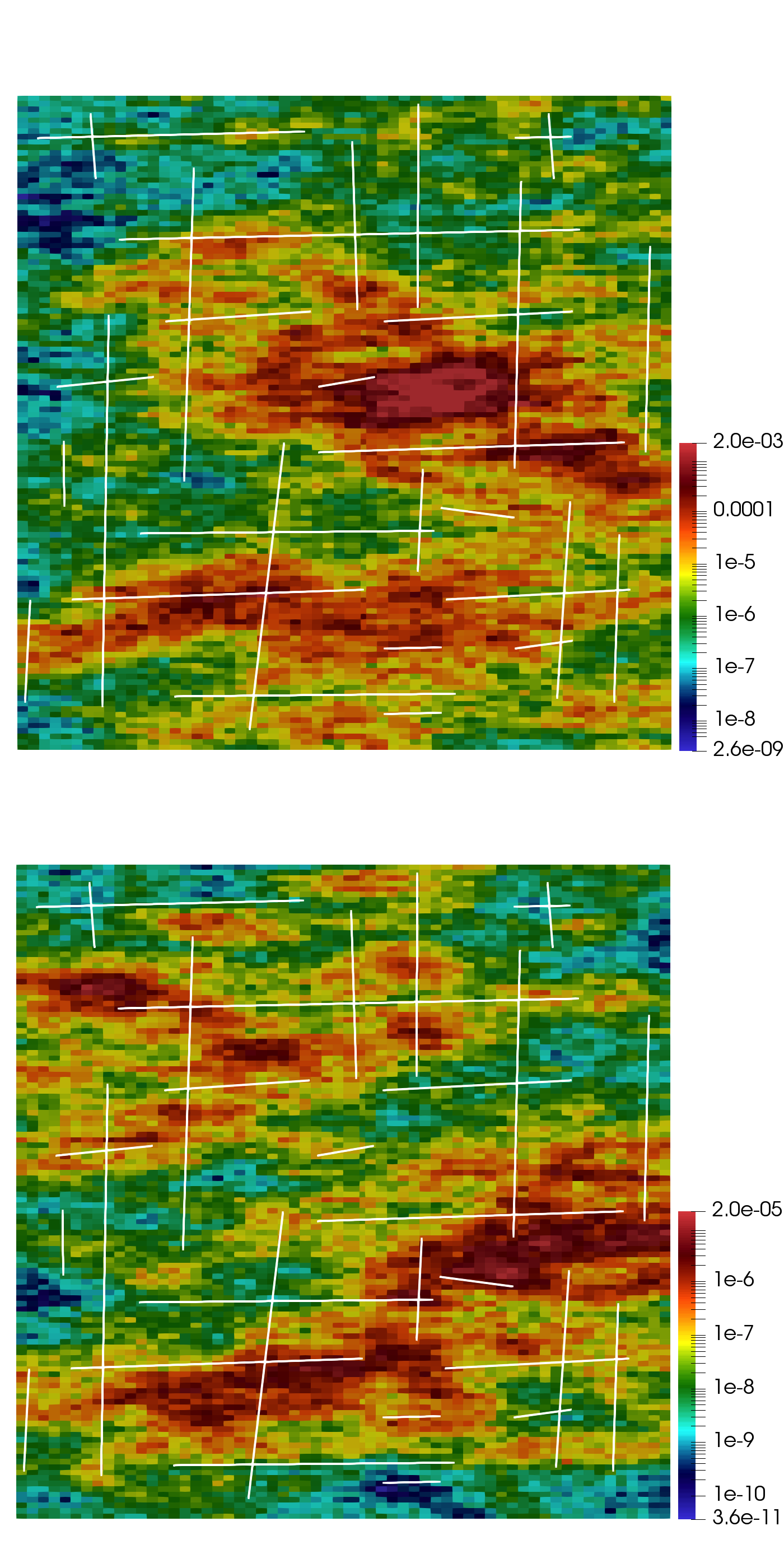}
\includegraphics[width=0.32 \textwidth]{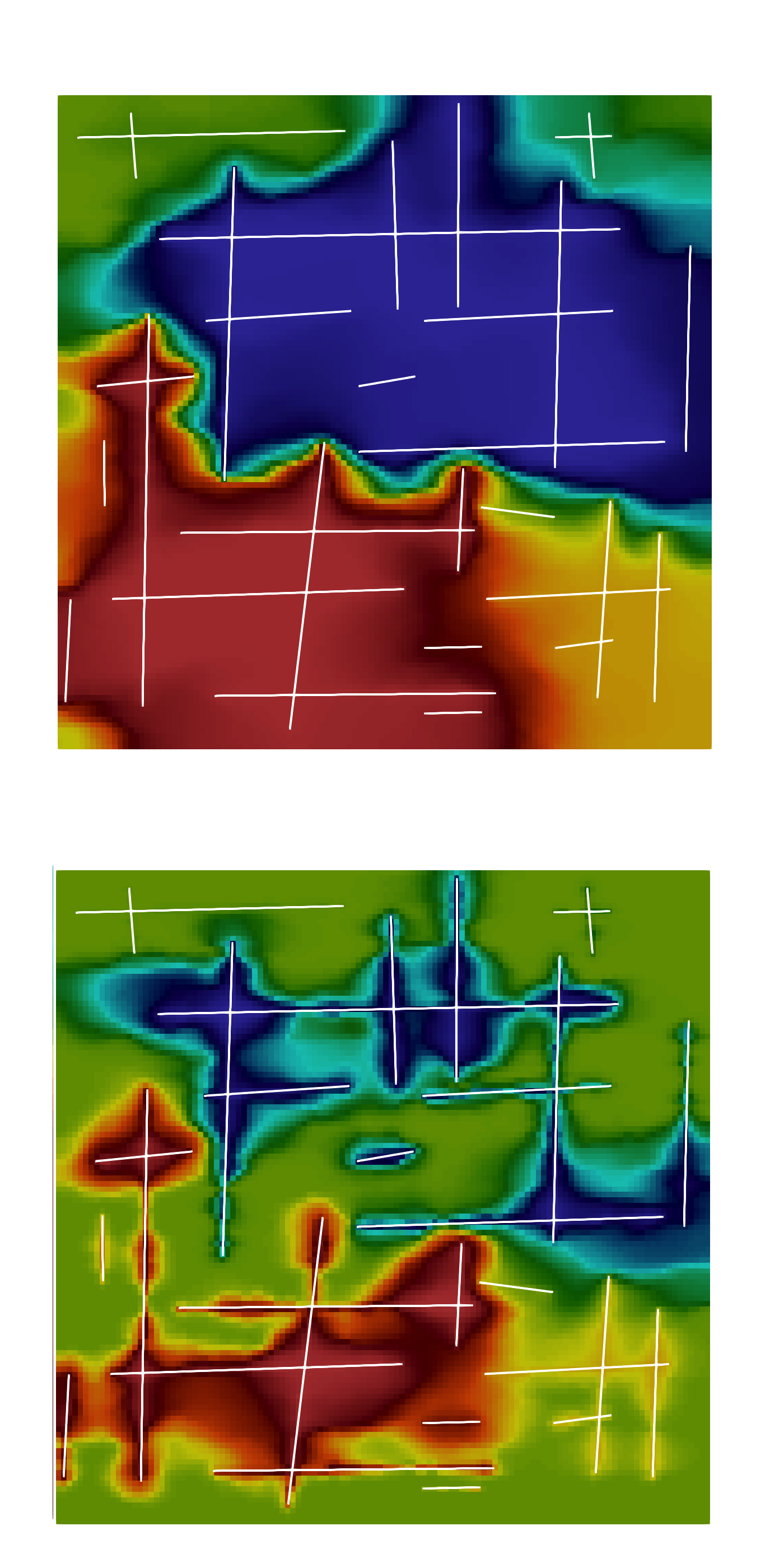}
\includegraphics[width=0.32 \textwidth]{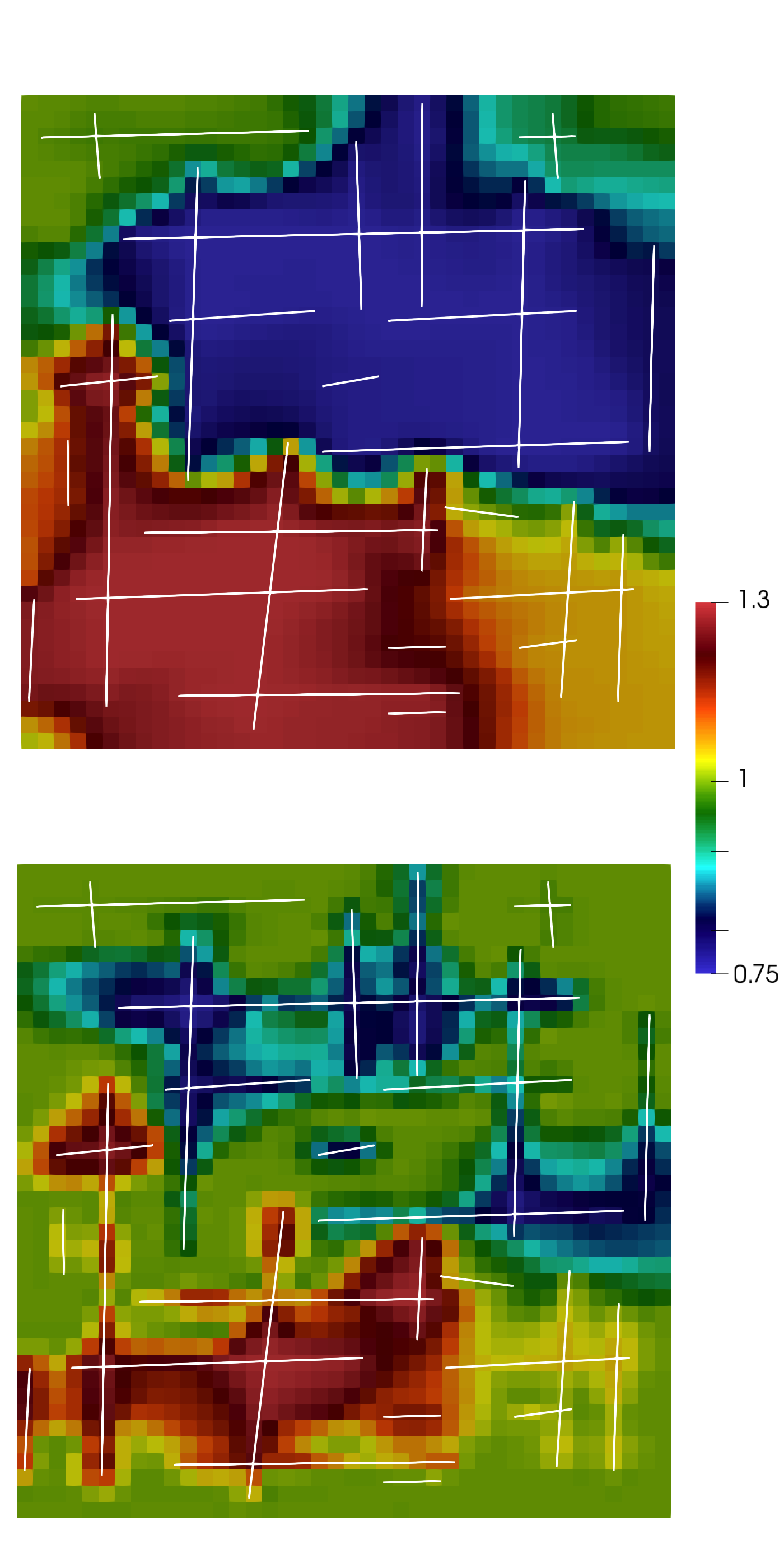}
\caption{Numerical results for pressure at final time for \textit{Test 2}.
Top: first continuum. Bottom: second continuum.
First column: heterogeneous permeabilities, $k_1$ and $k_2$.
Second column: fine-scale solution, $DOF_f = 29 842$.
Third column: coarse grid solution using NLMC, $DOF_c = 3 565$.
}
\label{fig:u2}
\end{figure}

\begin{table}[h!]
\centering
\begin{tabular}{ |c | c | c | c | c | }
\hline
$K^s$ &
$s = 1$&
$s = 2$&
$s = 3$&
$s = 4$ \\
\hline
\multicolumn{5}{|c|}{Coarse grid: $40 \times 40$} \\
\hline
$m = 5$ 		& 0.299 & 0.049 &	0.043 & 0.045	  \\ \hline
$m = 15$ 	& 0.821 & 0.165 &	0.051 & 0.061	  \\ \hline
$m = 25$ 	& 1.275 & 0.300 &	0.047 & 0.060	  \\ \hline
$m = 35$ 	& 1.641 & 0.422 &	0.047 & 0.054	  \\ \hline
$m = 50$ 	& 2.022 & 0.559 &	0.057 & 0.044	  \\ \hline
\multicolumn{5}{|c|}{Reconstructed fine grid} \\
\hline
$m = 5$ 		&  24.415 & 6.014 & 1.444 & 0.344	  \\ \hline
$m = 15$ 	&  24.292 & 5.999 & 1.446 & 0.344	  \\ \hline
$m = 25$ 	&  24.224 &	 5.996 & 1.447 & 0.341	  \\ \hline
$m = 35$ 	&  24.182 & 5.999	& 1.447 & 0.339	  \\ \hline
$m = 50$ 	&  24.145 &	 6.006 & 1.447 & 0.338	  \\ \hline
\end{tabular}
\,\,\,
\begin{tabular}{ |c | c | c | c | c | }
\hline
$K^s$ &
$s = 1$&
$s = 2$&
$s = 3$&
$s = 4$ \\
\hline
\multicolumn{5}{|c|}{Coarse grid: $40 \times 40$} \\
\hline
$m = 5$ 		& 0.195 & 0.104 &	0.100 & 0.099	  \\ \hline
$m = 15$ 	& 0.494 & 0.201 &	0.182 & 0.181	  \\ \hline
$m = 25$ 	& 0.746 & 0.270 &	0.230 & 0.229	  \\ \hline
$m = 35$ 	& 0.960 & 0.326 &	0.259 & 0.258	  \\ \hline
$m = 50$ 	& 1.210 & 0.390 &	0.277 & 0.276	  \\ \hline
\multicolumn{5}{|c|}{Reconstructed fine grid} \\
\hline
$m = 5$   	& 25.416 &	6.314 &	 1.601 & 0.577	  \\ \hline
$m = 15$ 	& 25.406 &	6.332 &	 1.665 & 0.716	  \\ \hline
$m = 25$ 	& 25.380 &	6.339 &	 1.677 & 0.741	  \\ \hline
$m = 35$ 	& 25.354 &	6.344 &	 1.673 & 0.733	  \\ \hline
$m = 50$ 	& 25.328 &	6.351 &	 1.661 & 0.713	  \\ \hline
\end{tabular}
\caption{Relative errors on the coarse grid and reconstructed fine-grid solutions for \textit{Test 2}. Coarse mesh: $40 \times 40$. $DOF_c  = 3565$ and $DOF_f = 29842$.
Left: first continuum. Right: second continuum. }
\label{err2}
\end{table}

The other model parameters used are as follows:
\begin{itemize}
\item \textit{Test 1}. Homogeneous matrix continuum properties with $k_1 = 0.5 \cdot 10^{-6}$, $k_2 = 10^{-5}$ and $k_f = 1$.
\item \textit{Test 2}. Heterogeneous matrix continuum properties (see Figure \ref{fig:u2}).
\end{itemize}
and $c_1 = 10^{-5}$, $c_2 = 10^{-5}$, $c_f = 10^{-6}$.
For mass transfer term, we set $\sigma_{12} = \sigma_{21} = k_1$, $\sigma_{1f} = \sigma_{f1} = k_1$ and $\sigma_{2f} = \sigma_{f2} = k_2$.
We set a point source on the fractures' continuum ($f$) at the two coarse cells $[0.15, 0.10] \times [0.15, 0.1]$ and $[0.65, 0.6] \times [0.9, 0.85]$ with $q = \pm10^{-3}$. As initial condition, we set $p_0 = 1$.
We simulate $t_{max} = 0.1$ with 50 time steps for upscaled and fine-scale solvers.

We depict a multiscale basis functions for \textit{Test 1} and \textit{Test 2} on the Figures \ref{fig:sch2} and  \ref{fig:sch3}, respectively.

To compare the results, we use the relative $L^2$ errors between fine grid in upscaled coarse grid models.
We calculate errors on coarse grid ($e^C$) and on fine grid ($e^F$)
\[
e^C = \frac{ || p_C - \bar{p} ||_{L^2} }{ || p_C ||_{L^2} },  \quad
e^F = \frac{ || p - \bar{p}_F ||_{L^2}  }{ || p ||_{L^2} },
\]
where $\bar{p}$ is the upscaled coarse grid solution, $\bar{p}_F = R^T \bar{p}$ is the downscaled of fine grid solution for $\bar{p}$, $p$ is the reference fine grid solution, $p_C$ is the coarse grid cell average for reference fine grid solution $p$ and
\[
|| p_C - \bar{p} ||^2_{L^2} =  \sum_K ( p^K_C - \bar{p}^K )^2,
\quad  p^K_C = \frac{1}{|K|} \int_K p \, dx.
\]
We calculate errors for background first and second continuums.

The fine-scale system has dimension $DOF_f = 29842$. The upscaled model has dimension $DOF_c = 993$ for coarse mesh with 400 cells ($20 \times 20$, and $DOF_c = 3565$ for coarse mesh with 1600 cells ($40 \times 40$).
In Figures \ref{fig:u1} and  \ref{fig:u2}, we present the fine scale solution for \textit{Test 1} and \textit{Test 2} with homogeneous and heterogeneous  matrix continuum properties at final time of simulation. The first row present solutions on the fine grid and on the second row we depicted a coarse grid upscaled solution using proposed method.  For basis calculations, we use oversampled domain $K^+$ with 4 coarse cells layers oversampling.
We observe  good accuracy of the proposed method with less than one percent of error for all test cases.

In Tables \ref{err1a} and \ref{err1}, we present relative errors for two coarse grids and for different numbers of oversampling layers $K^s$ with $s = 1,2,3$ and $4$ for \textit{Test 1} with homogeneous matrix properties.
The relative errors for \textit{Test 2} with heterogeneous properties, we present on Table \ref{err2} for different numbers of oversampling layers $K^s$ with $s = 1,2,3$ and $4$ on $40 \times 40$ coarse grid. We present relative errors for different time layers, $t_m$ with $m = 5, 15, 25, 35$ and $50$.
We notice a huge reduction of the system dimension and very small errors for unsteady mixed dimensional coupled system when we take sufficient number of the oversampling layers. For example, when we take only one oversampling layer for construction of the local domains, we obtain large errors for reconstructed coarse grid for both test cases. When we construct multiscale basis functions in $K^4$, we obtain less then one percent of errors for both continuum for the upscaled coarse grid solution and for reconstructed fine grid solution.

The fine-scale systems have size $DOF_f = 29842$ with solution time 69.74 seconds for \textit{Test 2}. The upscaled model  on $40 \times 40 $ coarse mesh has $DOF_c = 3565$ with solution time 4.84 seconds for $K^4$ and 3.85 seconds for $K^1$. We note that, the number of oversampling layers affect to the sparsity of the system but number of degrees of freedom is similar for any number of oversampling layers.  The proposed method is shown to be very efficient and provides good accuracy.


\section{Acknowledgements}
MV's  work is supported by the grant of the Russian Scientific Found N17-71-20055.
GP's work is supported by the mega-grant of the Russian Federation Government (N 14.Y26.31.0013).
EC's work is partially supported by Hong Kong RGC General Research Fund (Project 14304217) and CUHK Direct Grant for Research 2017-18.

\bibliographystyle{plain}
\bibliography{lit}

\begin{thebibliography}{10}

\bibitem{akkutlu2017multiscale}
I~Yucel Akkutlu, Yalchin Efendiev, Maria Vasilyeva, and Yuhe Wang.
\newblock Multiscale model reduction for shale gas transport in a coupled
  discrete fracture and dual-continuum porous media.
\newblock {\em Journal of Natural Gas Science and Engineering}, 2017.

\bibitem{akkutlu2018multiscale}
I~Yucel Akkutlu, Yalchin Efendiev, Maria Vasilyeva, and Yuhe Wang.
\newblock Multiscale model reduction for shale gas transport in poroelastic
  fractured media.
\newblock {\em Journal of Computational Physics}, 353:356--376, 2018.

\bibitem{akkutlu2012multiscale}
I~Yucel Akkutlu, Ebrahim Fathi, et~al.
\newblock Multiscale gas transport in shales with local kerogen
  heterogeneities.
\newblock {\em SPE journal}, 17(04):1--002, 2012.

\bibitem{akkutlu2015multiscale}
IY~Akkutlu, Yalchin Efendiev, and Maria Vasilyeva.
\newblock Multiscale model reduction for shale gas transport in fractured
  media.
\newblock {\em Computational Geosciences}, pages 1--21, 2015.

\bibitem{barenblatt1960basic}
GI~Barenblatt, Iu~P Zheltov, and IN~Kochina.
\newblock Basic concepts in the theory of seepage of homogeneous liquids in
  fissured rocks [strata].
\newblock {\em Journal of applied mathematics and mechanics}, 24(5):1286--1303,
  1960.

\bibitem{bosma2017multiscale}
Sebastian Bosma, Hadi Hajibeygi, Matei Tene, and Hamdi~A Tchelepi.
\newblock Multiscale finite volume method for discrete fracture modeling on
  unstructured grids (ms-dfm).
\newblock {\em Journal of Computational Physics}, 2017.

\bibitem{CYH2016adaptive}
Eric Chung, Yalchin Efendiev, and Thomas~Y Hou.
\newblock Adaptive multiscale model reduction with generalized multiscale
  finite element methods.
\newblock {\em Journal of Computational Physics}, 320:69--95, 2016.

\bibitem{chung2017coupling}
Eric~T Chung, Yalchin Efendiev, Tat Leung, and Maria Vasilyeva.
\newblock Coupling of multiscale and multi-continuum approaches.
\newblock {\em GEM-International Journal on Geomathematics}, 8(1):9--41, 2017.

\bibitem{chung2018nonfrac}
Eric~T Chung, Yalchin Efendiev, Wing~Tat Leung, Maria Vasilyeva, and Yating
  Wang.
\newblock Non-local multi-continua upscaling for flows in heterogeneous
  fractured media.
\newblock {\em Journal of Computational Physics}, 2018.

\bibitem{douglas1990dual}
Jim Douglas~Jr and T~Arbogast.
\newblock Dual porosity models for flow in naturally fractured reservoirs.
\newblock {\em Dynamics of Fluids in Hierarchical Porous Media}, pages
  177--221, 1990.

\bibitem{egh12}
Y.~Efendiev, J.~Galvis, and T.~Hou.
\newblock Generalized multiscale finite element methods.
\newblock {\em Journal of Computational Physics}, 251:116--135, 2013.

\bibitem{eh09}
Y.~Efendiev and T.~Hou.
\newblock {\em {Multiscale Finite Element Methods: Theory and Applications}},
  volume~4 of {\em Surveys and Tutorials in the Applied Mathematical Sciences}.
\newblock Springer, New York, 2009.

\bibitem{formaggia2014reduced}
Luca Formaggia, Alessio Fumagalli, Anna Scotti, and Paolo Ruffo.
\newblock A reduced model for darcy’s problem in networks of fractures.
\newblock {\em ESAIM: Mathematical Modelling and Numerical Analysis},
  48(4):1089--1116, 2014.

\bibitem{hkj12}
H.~Hajibeygi, D.~Kavounis, and P.~Jenny.
\newblock A hierarchical fracture model for the iterative multiscale finite
  volume method.
\newblock {\em Journal of Computational Physics}, 230(24):8729--8743, 2011.

\bibitem{houwu97}
T.~Hou and X.H. Wu.
\newblock A multiscale finite element method for elliptic problems in composite
  materials and porous media.
\newblock {\em J. Comput. Phys.}, 134:169--189, 1997.

\bibitem{jenny2005adaptive}
Patrick Jenny, Seong~H Lee, and Hamdi~A Tchelepi.
\newblock Adaptive multiscale finite-volume method for multiphase flow and
  transport in porous media.
\newblock {\em Multiscale Modeling \& Simulation}, 3(1):50--64, 2005.

\bibitem{karimi2003efficient}
Mohammad Karimi-Fard, Luis~J Durlofsky, Khalid Aziz, et~al.
\newblock An efficient discrete fracture model applicable for general purpose
  reservoir simulators.
\newblock In {\em SPE Reservoir Simulation Symposium}. Society of Petroleum
  Engineers, 2003.

\bibitem{karimi2001numerical}
Mohammad Karimi-Fard, Abbas Firoozabadi, et~al.
\newblock Numerical simulation of water injection in 2d fractured media using
  discrete-fracture model.
\newblock In {\em SPE annual technical conference and exhibition}. Society of
  Petroleum Engineers, 2001.

\bibitem{li2018multiscale}
Qiuqi Li, Yuhe Wang, and Maria Vasilyeva.
\newblock Multiscale model reduction for fluid infiltration simulation through
  dual-continuum porous media with localized uncertainties.
\newblock {\em Journal of Computational and Applied Mathematics}, 336:127--146,
  2018.

\bibitem{logg2009efficient}
Anders Logg.
\newblock Efficient representation of computational meshes.
\newblock {\em International Journal of Computational Science and Engineering},
  4(4):283--295, 2009.

\bibitem{logg2012automated}
Anders Logg, Kent-Andre Mardal, and Garth Wells.
\newblock {\em Automated solution of differential equations by the finite
  element method: The FEniCS book}, volume~84.
\newblock Springer Science \& Business Media, 2012.

\bibitem{lunati2006multiscale}
Ivan Lunati and Patrick Jenny.
\newblock Multiscale finite-volume method for compressible multiphase flow in
  porous media.
\newblock {\em Journal of Computational Physics}, 216(2):616--636, 2006.

\bibitem{martin2005modeling}
Vincent Martin, J{\'e}r{\^o}me Jaffr{\'e}, and Jean~E Roberts.
\newblock Modeling fractures and barriers as interfaces for flow in porous
  media.
\newblock {\em SIAM Journal on Scientific Computing}, 26(5):1667--1691, 2005.

\bibitem{geothermal2018multiscale}
Timothy Praditia, Rainer Helmig, and Hadi Hajibeygi.
\newblock Multiscale formulation for coupled flow-heat equations arising from
  single-phase flow in fractured geothermal reservoirs.
\newblock {\em Computational Geosciences}, pages 1--18, 2018.

\bibitem{sanchez1980non}
Enrique S{\'a}nchez-Palencia.
\newblock {Non-homogeneous media and vibration theory}.
\newblock In {\em {Non-homogeneous media and vibration theory}}, volume 127,
  1980.

\bibitem{talonov2016numerical}
Alexey Talonov and Maria Vasilyeva.
\newblock On numerical homogenization of shale gas transport.
\newblock {\em Journal of Computational and Applied Mathematics}, 301:44--52,
  2016.

\bibitem{ctene2016algebraic}
Matei {\c{T}}ene, Mohammed~Saad Al~Kobaisi, and Hadi Hajibeygi.
\newblock Algebraic multiscale method for flow in heterogeneous porous media
  with embedded discrete fractures (f-ams).
\newblock {\em Journal of Computational Physics}, 321:819--845, 2016.

\bibitem{ctene2017projection}
Matei {\c{T}}ene, Sebastian~BM Bosma, Mohammed~Saad Al~Kobaisi, and Hadi
  Hajibeygi.
\newblock Projection-based embedded discrete fracture model (pedfm).
\newblock {\em Advances in Water Resources}, 105:205--216, 2017.

\bibitem{vasilyeva2018constrained}
Maria Vasilyeva, Eric~T Chung, Yalchin Efendiev, and Jihoon Kim.
\newblock Constrained energy minimization based upscaling for coupled flow and
  mechanics.
\newblock {\em arXiv preprint arXiv:1805.09382}, 2018.

\bibitem{vasilyeva2018nonlocal}
Maria Vasilyeva, Eric~T Chung, Wing~Tat Leung, and Valentin Alekseev.
\newblock Nonlocal multicontinuum (nlmc) upscaling of mixed dimensional coupled
  flow problem for embedded and discrete fracture models.
\newblock {\em arXiv preprint arXiv:1805.09407}, 2018.

\bibitem{vasilyeva2018upscaling}
Maria Vasilyeva, Eric~T Chung, Wing~Tat Leung, Yating Wang, and Denis
  Spiridonov.
\newblock Upscaling method for problems in perforated domains with
  non-homogeneous boundary conditions on perforations using non-local
  multi-continuum method (nlmc).
\newblock {\em arXiv preprint arXiv:1805.09420}, 2018.

\bibitem{warren1963behavior}
JE~Warren, P~Jj Root, et~al.
\newblock The behavior of naturally fractured reservoirs.
\newblock {\em Society of Petroleum Engineers Journal}, 3(03):245--255, 1963.

\bibitem{wu2011multiple}
Yu-Shu Wu, Yuan Di, Zhijiang Kang, and Perapon Fakcharoenphol.
\newblock A multiple-continuum model for simulating single-phase and multiphase
  flow in naturally fractured vuggy reservoirs.
\newblock {\em Journal of Petroleum Science and Engineering}, 78(1):13--22,
  2011.

\bibitem{wu2007triple}
Yu-Shu Wu, Christine Ehlig-Economides, Guan Qin, Zhijang Kang, Wangming Zhang,
  Babatunde Ajayi, and Qingfeng Tao.
\newblock A triple-continuum pressure-transient model for a naturally fractured
  vuggy reservoir.
\newblock 2007.

\bibitem{xu2001modeling}
Tianfu Xu and Karsten Pruess.
\newblock Modeling multiphase non-isothermal fluid flow and reactive
  geochemical transport in variably saturated fractured rocks: 1. methodology.
\newblock {\em American Journal of Science}, 301(1):16--33, 2001.

\bibitem{yao2010discrete}
Jun Yao, Zhaoqin Huang, Yajun Li, Chenchen Wang, Xinrui Lv, et~al.
\newblock Discrete fracture-vug network model for modeling fluid flow in
  fractured vuggy porous media.
\newblock In {\em International oil and gas conference and exhibition in
  China}. Society of Petroleum Engineers, 2010.

\end{thebibliography}

\end{document}